\newif
\ifconfver
\confvertrue        %declaring conference version true
\ifconfver
    \documentclass[10pt,journal,final,twoside]{IEEEtran}
\else
    \documentclass[11pt,draftcls,onecolumn]{IEEEtran}
\fi
\usepackage{xspace,empheq,fancybox,amsmath,amssymb,epsfig,subfigure,syntonly,times,amsthm} \usepackage{psfrag,color,bm,array,cite}
\usepackage{url,cite,footnote,xspace,syntonly,algorithm,algorithmic}
\usepackage{verbatim,multirow,slashbox}
\usepackage[T1]{fontenc}

\newtheorem{remark}{\bfseries Remark}
\input{mysymbol.sty}

\newcommand{\diag}{\mathrm{diag}}

\newcommand{\hao}{\color{black}{}}
\newcommand{\ml}{\color{black}{}}
\begin{document}

%%%%%%%%%%%%%%%%%%%%%%%%%%%%%%%%%%%%%%%%%%%%%%%%%%%%%%%%%%%%%%%%%%%%%%
%                                                                    %
%               Paper Title                                          %
%                                                                    %
%%%%%%%%%%%%%%%%%%%%%%%%%%%%%%%%%%%%%%%%%%%%%%%%%%%%%%%%%%%%%%%%%%%%%%

\title{Fast Local Voltage Control under Limited Reactive Power: Optimality and Stability Analysis}
\author{
\IEEEauthorblockN{Hao Zhu},~\IEEEmembership{Member, IEEE}, and \IEEEauthorblockN{Hao Jan Liu},~\IEEEmembership{Student Member, IEEE}

\thanks{\protect\rule{0pt}{3mm} This work is partially supported by the Illinois Center for a Smarter Electric Grid (ICSEG) and Trustworthy Cyber Infrastructure for the Power Grid (TCIPG).} \thanks{\protect\rule{0pt}{3mm} The authors are with the Department of Electrical \& Computer Engineering, University of Illinois, 306 N. Wright Street, Urbana, IL, 61801, USA; Emails: {\{haozhu,haoliu6\}{@}illinois.edu}.
}}

% % The paper headers
\markboth{IEEE TRANSACTIONS ON POWER SYSTEMS (SUBMITTED)}%
{Zhu \MakeLowercase{\textit{et al.}}: Fast Local Voltage Control Under Limited Reactive Power}
\renewcommand{\thepage}{}
\maketitle
\pagenumbering{arabic}

%%%%%%%%%%%%%%%%%%%%%%%%%%%%%%%%%%%%%%%%%%%%%%%%%%%%%%%%%%%%%%%%%%%%%%
%                                                                    %
%                   Abstract                                         %
%                                                                    %
%%%%%%%%%%%%%%%%%%%%%%%%%%%%%%%%%%%%%%%%%%%%%%%%%%%%%%%%%%%%%%%%%%%%%%
%
\begin{abstract}
High penetration of distributed energy resources presents several challenges and opportunities for voltage regulation in power distribution systems. A local reactive power (VAR) control framework will be developed that can fast respond to voltage mismatch and address the robustness issues of (de-)centralized approaches against communication delays and noises. Using local bus voltage measurements, the proposed gradient-projection based schemes explicitly account for the VAR limit of every bus, and are proven convergent to a surrogate centralized problem with proper parameter choices. This optimality result quantifies the capability of local VAR control without requiring any real-time communications. The proposed framework and analysis generalize earlier results on the droop VAR control design, which may suffer from under-utilization of VAR resources in order to ensure stability. Numerical tests have demonstrated the validity of our analytical results and the effectiveness of proposed approaches implemented on realistic three-phase systems. 

\end{abstract}

%\begin{IEEEkeywords}
%\end{IEEEkeywords}

%\newpage

%%%%%%%%%%%%%%%%%%%%%%%%%%%%%%%%%%%%%%%%%%%%%%%%%%%%%%%%%%%%%%%%%%%%%%
% %%
% %%      Section: Intro %%
% %%%%%%%%%%%%%%%%%%%%%%%%%%%%%%%%%%%%%%%%%%%%%%%%%%%%%%%%%%%%%%%%%%%%
%%%

\section{Introduction}\label{sec:intro}

Recent technological advances in distribution systems have led to growing interest in the voltage regulation problem using distributed energy resources (DERs) of inverter interfaces \cite{Tur_proc11}. High variability of photovoltaic (PV) solar generation and abrupt load increase due to electric vehicle charging could result in unexpected voltage fluctuations, at  time-scales much faster than the existing voltage control using on-load tap-changers or capacitor banks; see e.g., \cite{Car_tps08,Tur_proc11,JahAli_tps13,Rob_tps13}. Further DER deployment gives rise to unprecedented capability of extremely fast voltage regulation through inverter VAR output control, in order to meet the 5$\%$ voltage deviation bounds per the ANSI Standard C84.1 \cite{ANSI}. 

With the full system information available centrally, the VAR control problem can be cast as an optimal power flow (OPF) one that minimizes system operational costs such as power losses and voltage violation \cite{Far_pesgm12}. Several distributed optimization algorithms have also been proposed to solve this centralized problem using information exchanges among neighboring buses. To balance overall system VAR resources, the second-stage VAR control in \cite{Rob_tps13} relies on a consensus averaging protocol. Alternating-direction method-of-multipliers (ADMM) has been advocated in \cite{dall_tsg13,lanl}, while a subgradient iterative solver has been developed by \cite{lam_tps15}. More recently, a stochastic-approximation approach has been adopted in \cite{vkgg_tps15} to handle high system variability and measurement noises. Nonetheless, all  (de-)centralized approaches would require high-quality communication of the measurement and control signals, which is not yet a reality for almost all distribution systems. Since these optimization-based control methods are designed in an open-loop fashion, potential communication delays or noises would challenge their optimality and stability for real-time implementations.
%would greatly challenge their implementations in real time, and even raise stability concerns since the control signals are optimized in an open-loop fashion. 

To tackle this, one can design VAR control strategies using locally available information such as bus voltage magnitude measurements \cite{Car_tps08,Tur_proc11,yeh_tps12}. Since power system voltage is {\ml more significantly affected} by local VAR inputs {\ml compared to those elsewhere} \cite[Sec. 10.8]{PSA_book}, a local control framework would very fast and effectively respond to voltage deviation. {\ml One challenge in local control design is the guaranteed dynamic system stability. As shown by \cite{qz_iet07,Kraiczy_PVcon,gvtvc_tsg13}, local inverter VAR control would potentially increase the number of tap changes triggered for on-load tap-changing (OLTC) transformers. Hence, stability of local VAR control could be crucial to ensure that it will not offset, or even adversely affect OLTC transformers and other voltage regulating devices.} 
One popular local control is the piecewise linear droop design using instantaneous local bus voltage mismatch input, as advocated by the IEEE 1547.8 Standard \cite{ieee1547}. However,  as shown by \cite{FarLow_cdc13} the droop slope has to be small enough to ensure system stability. This condition is equivalent to enforcing high penalty on inverter VAR output, which would result in insufficient utilization of VAR resources. A delayed droop scheme was developed in \cite{JahAli_tps13}, which could relax the droop stability condition. However, it is not clear how to choose the delay parameter to balance the stability and the convergence speed. Stability and optimality analysis for similar integral control approaches have been offered in \cite{bz_naps13,li_allerton14}, but they assume unconstrained VAR resources at every bus and thus fail to account for inverter rating limits in practice. 

 The present paper offers a general framework for developing local VAR control strategies that explicitly account for the VAR limits at every bus and potential VAR supplying penalty. This problem turns out to be a box-constrained quadratic optimization problem, which motivates us to leverage the iterative gradient-projection (GP) method \cite[Sec. 2.3]{Ber_NPbook}. Interestingly, the GP iterations naturally decouple into local updates requiring only bus voltage magnitude, shown to attain the optimal solution to a surrogate centralized VAR control problem under proper choice of parameters. Compared to the existing literature on (local) VAR control, our contribution is three-fold. First, our proposed GP-based local control generalizes the existing droop and delayed droop methods for limited VAR resources in \cite{FarLow_cdc13,JahAli_tps13}. Second, compared to these earlier approaches, our GP-based update can be stabilized under any choice of VAR supplying penalty. Last but for least, motivated by the acceleration from Newton's method, we have advocated to diagonally scale the GP stepsize parameters based on the inverse Hessian matrix. This proposed stepsize selection method can improve the system conditioning and convergence speed, while sacrificing no optimality condition. 

The rest of the paper is organized as follows. Sec. \ref{sec:ps} presents the linearized distribution system flow model, which is further reformulated by introducing graph matrices in Sec. \ref{sec:graph}. The gradient-projection (GP) algorithm is introduced in Sec. \ref{sec:gp} to solve the limit-constrained VAR optimization problem, which leads to equivalent local voltage control updates. 
%The fixed point of the GP iterations is shown to attain the optimal solution to a surrogate centralized problem, which points out the optimality trade-off of the control strategies that only require local voltage information. 
Stability conditions have been offered in Sec. \ref{sec:local} for selecting GP stepsize parameters, which are based on only system topology and line admittance. %Interestingly, the GP based framework can generalize earlier results for (delayed) droop control in \cite{FarLow_cdc13,JahAli_tps13}.  
Numerical tests on realistic (three-phase) feeder systems have been performed for real-time implementations, which demonstrate the effectiveness of the proposed scaled control design with improved stability and voltage regulation performance. {\ml An important issue of local VAR control is its interactions with other voltage regulation devices, and we are actively pursuing this direction to extend the present work.}

%%%%%%%%%%%%%%%%%%%%%%%%%%%%%%%%%%%%%%%%%%%%%%%%%%%%%%%%%%%%%%%%%%%%%%
% %%
% %%      Section: PS %%
% %%%%%%%%%%%%%%%%%%%%%%%%%%%%%%%%%%%%%%%%%%%%%%%%%%%%%%%%%%%%%%%%%%%%
%%%

\section{System Modeling and Problem Statement}\label{sec:ps}

\begin{figure}[tb]
\centering
\vspace{-3pt}
\includegraphics[width=.96\linewidth,clip = true, trim = .8in 2.88in  1.5in  2.36in]{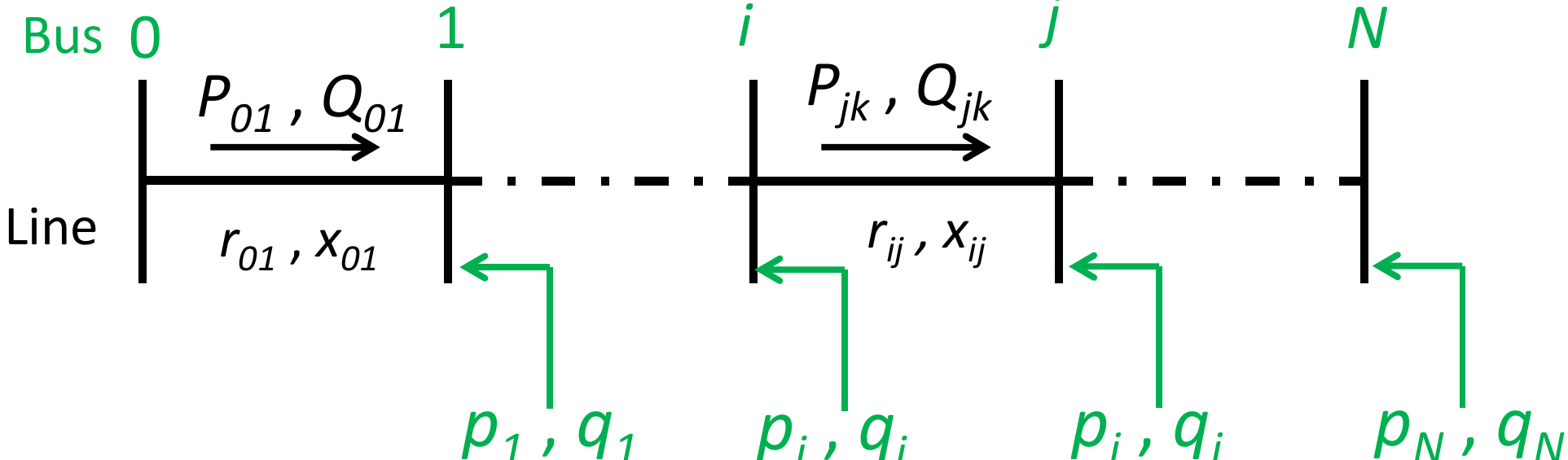}
\caption{A radial distribution feeder with bus and line associated variables.}
\label{fig:radial}
\vspace{-4mm}
\end{figure}

Consider a distribution network with  $N+1$ \textit{buses} collected in the set $\ccalN := \{0, 1,\ldots,N\}$, and \textit{line segments} represented by the set $\ccalL := \{(i,j)\} \subset \ccalN  \times \ccalN$; see Fig. \ref{fig:radial} for a radial feeder illustration. For tree-topology distribution networks, the number of line segments $|\ccalL| = (N+1)-1=N$. Bus $0$ denotes the point of common coupling (PCC), usually at the distribution substation and assumed to be of {\ml reference} voltage. For every bus $i$, let $V_i$ denote its voltage magnitude, and $p_i$ and $q_i$ denote the bus real and reactive power injection, all in per unit (p.u.). For each line segment $(i,j)$, let $r_{ij}$ and $x_{ij}$ denote its resistance and reactance, and $P_{ij}$ and $Q_{ij}$ the real and reactive power from bus $i$ to $j$, respectively. In addition, the subset $\ccalN_j \subset\ccalN$ denotes bus $j$'s neighboring buses that are further  down from the feeder head. The DistFlow equations \cite{mbfw_tpd89} to model the distribution network flow are given for every line $(i,j)\in \ccalL$ as 
\begin{subequations}
\label{distflow}
\begin{align}
P_{ij} \! -\! \sum_{k\in\ccalN_j} \!P_{jk} &= -p_j + r_{ij} \frac{P_{ij}^2+Q_{ij}^2}{V_i^2}, \label{dfP} \\
Q_{ij}\! -\! \sum_{k\in\ccalN_j} \!  Q_{jk} &= - q_j + x_{ij} \frac{P_{ij}^2+Q_{ij}^2}{V_i^2}, \label{dfQ} \\
V_i^2 -V_j^2 & =  2(r_{ij}P_{ij} + x_{ij}Q_{ij}) - (r_{ij}^2 +x_{ij}^2) \frac{P_{ij}^2+Q_{ij}^2}{V_i^2} \label{dfV}
\end{align}
\end{subequations}
where the nonlinear term $(P_{ij}^2+Q_{ij}^2)/V_i^2$ represents the squared line current magnitude, leading to the power loss terms in \eqref{distflow}.
Assuming the loss is negligible compared to line flow, a linear approximation of \eqref{distflow} can be constructed. The approximation error introduced is relatively small, at the order of 1\%  \cite{FarLow_cdc13}. Under relatively flat voltage profile, i.e., $V_i \approx 1,$ $\forall i \in \ccalN$, we have $V_i^2 -V_j^2 \approx 2(V_i -V_j)$. This leads to a small approximation error at about 0.25\% (1\%) if there is a 5\% (10\%) deviation in voltage magnitude approximation \cite{FarLow_cdc13}. Under the two assumptions,  the linearized DistFlow (LinDistFlow) model can be established for every $(i,j) \in \ccalE$, as 
\begin{subequations}
\label{lindistflow}
\begin{align}
\textstyle P_{ij} - \sum_{k\in\ccalN_j} P_{jk} &= -p_j, \label{ldfP} \\
\textstyle Q_{ij} - \sum_{k\in\ccalN_j} Q_{jk} &= -q_j, \label{ldfQ} \\
V_i -V_j & =  r_{ij}P_{ij} + x_{ij}Q_{ij}. \label{ldfV}
\end{align}
\end{subequations}

The total injected reactive power $q_j = q_j^g - q_j^c$, where $q_j^g$ denotes the reactive power {\ml contributed} by e.g., PV inverters per bus $j$,  while $q_j^c$ corresponds to the load reactive power consumption at bus $j$. Using the bold symbols to denote their vector counterparts, the VAR control problem becomes to solve for $\bbq^g$ for given $\bbp$ and $\bbq^c$, such that 
\begin{subequations}
\label{stat}
\begin{align}
\min_{\bbV, \bbq^g, \bbP, \bbQ}~~ &~~ \frac{1}{2}\|\bbV - \bbmu\|_2^2 \label{statf}\\
\mathrm{subject~to}~~&~~  \underline{q}_j\leq q_j^g \leq \bbarq_j, ~~ \forall j \label{statq}\\
&~~V_0 = 1,~\textrm{and}~\eqref{ldfP}-\eqref{ldfV} \label{statp}
\end{align}
\end{subequations}
where the substation voltage $V_0$ is fixed to be unit and $\bbmu$ is the voltage profile to be achieved. For example, to attain the flat-voltage profile, we can set $\bbmu=\mathbf 1$. Other choices also exist depending on certain operational objectives such as better energy-efficiency as motivated by conservation voltage reduction (CVR); see e.g., \cite{liu_sgcomm14}. The bounds on $q_j^g$ could be due to the inverters' apparent power limit and real power generation \cite{Tur_proc11}, or depend on the inverter power factor ratings. %which will make the VAR control problem more interesting. 

%Since constraints \eqref{statq}-\eqref{statp} are linear, \eqref{stat} is a convex quadratic program (QP). It can be easily solved with the knowledge on the full system model, as well as the real and reactive power load everywhere. Such a centralized framework however, has a high communication overhead, which is not yet realistic or widely available for existing feeders. Meanwhile, noise or delays of information exchange could lead to potential instability concerns and fail to cope with the increasing variability in generation and loads. To address these challenges, we will design control strategies that can quickly manage the reactive power based on local voltage measurements. 

%%%%%%%%%%%%%%%%%%%%%%%%%%%%%%%%%%%%%%%%%%%%%%%%%%%%%%%%%%%%%%%%%%%%%%
%% %%%
% %%%      Section: Method %%%
% %%%%%%%%%%%%%%%%%%%%%%%%%%%%%%%%%%%%%%%%%%%%%%%%%%%%%%%%%%%%%%%%%%%%
%%%%
\section{Graph-Representation based VAR Control}\label{sec:graph}

{\ml We will present our new LinDistFlow representation using graph-based matrices, to better formulate the VAR control problem.} Let matrix $\bbM^o$ of size $(N+1)\times N$ denote the graph incidence matrix  for $(\ccalN, \ccalL)$; see e.g., \cite[pg. 6]{west_gt}. Its $\ell$-th column corresponds to one line segment $(i,j)\in\ccalL$, the entries of which are all zero except for the $i$-th and $j$-th ones. Specifically, set $M^o_{il} =1$ and $M^o_{jl} = -1$ when $j\in\ccalN_i$; i.e., bus $i$ is closer to the feeder head. Furthermore, let vector $\bbm_0^T$ denote the first row of $\bbM^o$ which corresponds to bus 0, while the rest of matrix denoted by $\bbM$ of size $N\times N$; i.e., $\bbM^o = [\bbm_0~\bbM^T]^T$. Since the tree topology of $(\ccalN, \ccalL)$ ensures that it is a connected graph, the rank of $\bbM^o$ equal to $(N+1)-1$ = $N$ \cite{west_gt}. Therefore, the square matrix $\bbM$ is of full rank $N$, and thus invertible.

With the graph matrix notation and $V_0 =1$, all voltage equations \eqref{ldfV} in the LinDistFlow model become
\begin{align}
(\bbM^o)^T [V_0~\bbV^T]^T = \bbm_0 + \bbM^T \bbV  &=  \bbD_r \bbP + \bbD_x \bbQ \label{ldfVm}
\end{align}
where $\bbD_r$ is an $N\times N$ diagonal matrix with the $\ell$-th diagonal entry equal to $r_{ij}$; and similarly for $\bbD_r$ which captures all $x_{ij}$'s. Similarly, the power balance equations in \eqref{ldfP} and \eqref{ldfQ} can be respectively concatenated into
\begin{align}
- \bbM \bbP = - \bbp, \label{ldfPm}\\
-\bbM \bbQ = - \bbq. \label{ldfQm}
\end{align}
Solving for $\bbP$ and $\bbQ$ and substituting them into \eqref{ldfVm} yield
\begin{align}
\bbM^T \bbV = \bbD_r \bbM^{-1} \bbp +\bbD_x \bbM^{-1} \bbq - \bbm_0, \label{ldfVm1}
\end{align}
or equivalently,
\begin{align}
 \bbV = \bbR \bbp +\bbX \bbq - \bbM^{-T}\bbm_0 \label{ldfVm2}
\end{align}
where the two invertible matrices $\bbR:=\bbM^{-T}\bbD_r \bbM^{-1}$ and $\bbX := \bbM^{-T}\bbD_x \bbM^{-1}$. Since $\bbq = \bbq^g -\bbq^c$, the linear relation from input $\bbq^g$ to the output $\bbV$ becomes
\begin{align}
\bbV =  \bbR \bbp +\bbX \bbq^g  - \bbX \bbq^c - \bbM^{-T}\bbm_0  = \bbX \bbq^g + \bar{\bbV}  \label{ldfVm3}
\end{align}
where $\bar{\bbV} : = \bbR \bbp - \bbX \bbq^c - \bbM^{-T}\bbm_0$ denotes the voltage profile under no additional VAR support.  

\begin{proposition}
Both $\bbR$ and $\bbX$ are positive definite (PD).
\end{proposition}

{\hao {\it Proof:} Picking any non-zero vector $\bbz$ of length $N$ and defining $\bbz' := \bbM^{-1}\bbz$, one can show that  
\begin{align*}
\ml{\bbz^T\bbR\bbz = (\bbM^{-1}\bbz)^T \bbD_r (\bbM^{-1}\bbz) = (\bbz')^T \bbD_r (\bbz').}
\end{align*}
Since diagonals of $\bbD_r$ are all positive resistance values, $\bbz^T\bbR\bbz>0$ holds for any non-zero $\bbz$. This way, we prove $\bbR$ is a PD matrix. Similarly, since all diagonals of $\bbD_x$ are positive, matrix $\bbX$ is PD as well. \qed}

Note that the matrix LinDistFlow model \eqref{ldfVm2} has been established by \cite{FarLow_cdc13}. However, the linearity and positive definiteness were shown using a more complicated induction-based proof. Using the graph incidence matrix, the analysis here for the matrix LinDistFlow model is much more straightforward. Moreover, the definitions based on $\bbM$ suggest a very efficient way to form $\bbR$ and $\bbX$, and even to model the approximation error introduced by perturbed system information. Interestingly, the inverse of $\bbX$ as denoted by $\bbB : = \bbX^{-1} = \bbM \bbD_x^{-1}\bbM^T$ is also a PD matrix, and it is actually the power network \textit{Bbus matrix} used in  the fast decoupled power flow (FDPF) analysis; see e.g., \cite[Sec. 6.16]{Wollenberg-book}. The matrix LinDistFlow model \eqref{ldfVm3} is equivalent to
\begin{align}
\bbq^g  = \bbX^{-1} (\bbV- \bar{\bbV}) = \bbB  (\bbV- \bar{\bbV}).  \label{ldfVm4}
\end{align}

\begin{remark} \label{rmk:loss} (Power loss terms.) {\rm 
Albeit an approximation of \eqref{distflow}, the linear relation in \eqref{ldfVm4} is very meaningful as it could be considered as the sensitivity of $\bbq^g$ to $\bbV$ around any operating point. As shown by \cite{bz_naps13,lanl}, change of reactive power injection does not affect very much the power loss terms in \eqref{dfP}-\eqref{dfV}. Hence, relatively constant loss terms can be even captured by $\bbarbbV$ as an operating-point related voltage profile. }
\end{remark}

\begin{remark} \label{rmk:vol2} (Squared voltage profile.) {\rm The other assumption  used for linearization  relates to the squared voltage terms. However, \eqref{ldfVm4} can be generalized to include the original squared voltage, instead of $\bbV$, as the input. To keep the squared voltage difference $V_i^2-V_j^2$ of \eqref{dfV},  the graph-incidence based reformulation \eqref{ldfV} and \eqref{ldfVm} can be modified by substituting the voltage squared vector, with a scaling difference by a factor of $1/2$. This way, matrix $\bbX$ becomes the linear sensitivity of (squared) voltage profile due to change in $\bbq^g$. The exact form of this sensitivity matrix is given by \cite{bz_naps13}, which is closely approximated by $\bbX$. Based on this modification, the VAR control problem \eqref{stat} and the ensuing analysis can be generalized accordingly to use the squared voltage term for higher accuracy.}
\end{remark}

\begin{remark} \label{rmk:meshed} {\hao (Meshed networks.) {\rm Although the matrix LinDistFlow model in \eqref{ldfVm4} has been derived for tree-topology networks, it can also be generalized to meshed networks. Since $\bbB:=\bbM \bbD_x^{-1}\bbM^T$, it is the weighted graph Laplacian matrix. This exactly coincides with the definition of the system Bbus matrix in the popular fast-decoupled power flow (FDPF) model for transmission network analysis \cite[Sec. 6.16]{Wollenberg-book}. Hence, the model \eqref{ldfVm4} and the ensuing algorithms also hold for general distribution networks such as ring-topology systems.}}
\end{remark}

Remarks \ref{rmk:loss}-\ref{rmk:meshed} corroborate the validity of using the linearized model \eqref{ldfVm4} for representing realistic distribution networks with lossy lines, non-flat voltage profile, and even meshed topology. This will be further demonstrated by numerical tests for even three-phase unbalanced cases in Sec. \ref{sec:sim}.

Upon defining the voltage mismatch $\tdbbV:= \bbmu-\bar{\bbV}$ and representing $\bbq^g$ by $\bbq$ for notational convenience in the rest of the paper, the VAR control problem now becomes
\begin{subequations}
\label{cenq}
\begin{align}
\bbq^\dagger = \arg \min_{\bbq}~~ &~ \frac{1}{2}\|\bbV - \bbmu\|^2= \frac{1}{2}\|\bbX \bbq - \tdbbV\|^2 \label{cenqf}\\
\mathrm{subject~to}~~&~~  \underline{\bbq}\leq \bbq \leq \bar{\bbq} \label{cenqq}
\end{align}
\end{subequations}
where $\|\cdot\|$ denotes the Euclidean norm operator. 
This is a box-constrained  quadratic program (QP) and could be easily solved using off-the-shelf convex solvers.  %From this point on, $\bbq^g$ will be just denoted by $\bbq$ with no further notational confusion. 
To account for the cost of supplying VAR, a quadratic penalty term per bus will be introduced, as motivated by the linear droop control design \cite{FarLow_cdc13}. The total cost of supplying $\bbq$ is defined as
\begin{align}
C(\bbq) := %\sum_{j=1}^N C_j(q^g_j) = 
\sum_{j=1}^N c_j q_j^2 = \bbq^T \bbC \bbq \label{costqg}
\end{align}
where the penalty coefficient $c_j \geq 0$  and the diagonal matrix $\bbC :=\diag(c_1, \ldots, c_N)$. {\ml Clearly, the case of no penalty on VAR support as in \eqref{cenq} can be included by setting $\bbC = \bb0$. The reason for penalizing the VAR supply with a positive $c_j$ is three-fold. First, it will be shown in Sec. \ref{sec:droop} that the popular linear droop control design naturally comes from a positive $c_j$ which corresponds to the inverse of the droop slope. Second, since larger reactive power flow could result in higher line current level and hence power losses, it is possible to prevent unnecessary VAR loss and decrease the overall network power losses by discouraging the net VAR injection at some locations. Last but not least, as detailed soon in Sec. \ref{sec:local}, a larger $\bbC$ tends to better and faster stabilize the system. Thus, it would decrease the occurrence and level of abrupt changes in system dynamics, and has the potential to prevent local control actions from adversely affecting the operations of OLTC and other voltage regulating devices. }

To allow for local control schemes,  it turns out that the voltage mismatch norm in \eqref{cenq} needs to be weighted by the PD matrix $\bbB$, leading to the following surrogate VAR control problem 
%
%\begin{subequations}
\begin{align}
\bbq^\star = & \arg \min_{\underline{\bbq}\leq \bbq \leq \bar{\bbq}}~  f(\bbq) 
\label{cenqn}
\end{align}
where the objective 
\begin{align}
f(\bbq) := \frac{1}{2}(\bbX \bbq - \tdbbV)^T\bbB(\bbX \bbq - \tdbbV)  + \frac{1}{2} \bbq^T \bbC \bbq. \label{cenqfn}
\end{align}
This surrogate problem is still convex as $\bbB$ is PD. {\hao Because of the box constraint,  the weighted error norm with $\bbC = \bb0$ would attain a different solution compared to the original unweighted problem \eqref{cenq}. However, if every bus has unlimited VAR capability (i.e., unconstrained case), the optimal solutions to both error norm objectives coincide at $\bbq^\dag = \bbq^\star = \bbX^{-1}( \bbmu-\bar{\bbV})$  if $\bbC = \bb0$. This implies that under abundant VAR resources, the optimal solution to the weighted \eqref{cenqn} has the potential to closely approximate the minimum of the \eqref{cenq}. As detailed soon, the surrogate problem \eqref{cenqn} will facilitate the development of fast local control schemes.}

%%%%%%%%%%%%%%%%%%%%%%%%%%%%%%%%%%%%%%%%%%%%%%%%%%%%%%%%%%%%%%%%%%%%%%
%% %%%
% %%%      Section: Local %%%
% %%%%%%%%%%%%%%%%%%%%%%%%%%%%%%%%%%%%%%%%%%%%%%%%%%%%%%%%%%%%%%%%%%%%
%%%%
\section{Gradient-Projection Method} \label{sec:gp}
{\ml This section will introduce our proposed VAR control framework by solving the constrained optimization problem \eqref{cenqn}.}
The key of solving \eqref{cenqn} lies in the separable structure of the constraint. To project a vector to the set  $[\underline{\bbq}, \bar{\bbq}]$, one can threshold it on each coordinate. {\hao This motivates one to adopt the gradient-projection (GP) method, a generic optimization solver for constrained problems; see e.g., \cite[Sec. 2.3]{Ber_NPbook}. The GP method extends the iterative first-order gradient methods to solving constrained optimization problems like \eqref{cenqn}. Upon forming the gradient direction of \eqref{cenqfn}
by using $\bbB\bbX=\bbI$ as
\begin{align}
\nabla f(\bbq) := (\bbX \bbq - \tdbbV) +  \bbC \bbq, \label{graf}
\end{align}
the simplest GP iteration finds the feasible direction by projecting the gradient update, as
\begin{align}
\bbq(t+1) = \mathbb P \!\left[\bbq(t) - d\nabla f(\bbq(t)) \right]
%= \mathbb P [\bbq^g(t)-d (\bbX\bbq^g(t) - \tdbbV +  \bbC \bbq^g(t))]
\label{graprojs}
\end{align}
where the $\mathbb P$ operator thresholds any input to the constraint set $[\underline{\bbq}, \bar{\bbq}]$, and $d>0$ is the stepsize. Clearly, the GP iteration \eqref{graprojs} boils down to the steepest descent method if the optimization problem is unconstrained.

A more general GP form takes the weighted linear combination of the last iterate and the projection, given by
\begin{align}
\bbq(t+1)\! = \![1-\alpha(t)]\bbq(t) + \alpha(t)\mathbb P \!\left[\bbq(t) - d\nabla f(\bbq(t)) \right]
%= \mathbb P [\bbq^g(t)-d (\bbX\bbq^g(t) - \tdbbV +  \bbC \bbq^g(t))]
\label{graproj}
\end{align}
where the weighting parameter $\alpha(t)\in (0,1]$. This choice of $\alpha(t)$ ensures that $\bbq(t+1)$ would always satisfy the box constraint for any iteration $t$, as long as the last iterate $\bbq(t) \in [\underline{\bbq}, \bar{\bbq}]$ . Hence, every GP iterate is feasible to \eqref{cenqn} as long as $\bbq(0) \in [\underline{\bbq}, \bar{\bbq}]$. By scaling each bus $j$ with a different stepsize $d_j>0$, the most general GP form  is
\begin{align}
\bbq(t+1) \! = \! [1-\alpha(t)]\bbq(t) + \alpha(t)\mathbb P \!\left[\bbq(t) - \bbD\nabla f(\bbq(t)) \right]
\label{graprojd}
\end{align}
where the diagonal matrix $\bbD:=\diag(d_1,\ldots, d_N)$. Clearly, the original GP iteration \eqref{graproj} is a special case of \eqref{graprojd} by setting $\bbD = d\bbI$. 

Interestingly, the GP iteration \eqref{graprojd} can be easily implemented by setting the VAR input at all buses to be the instantaneous $\bbq(t)$. Since  $\bbq(t) \in [\underline{\bbq}, \bar{\bbq}]$ always holds , it is feasible to use the latest GP iterate as the network reactive power input by setting $\bbq^g = \bbq(t)$ at every iteration $t$. Under this setting,  the gradient direction for any $\bbq(t)$ becomes [cf. \eqref{ldfVm3}]
\begin{align}
\nabla f(\bbq(t))\! := \!(\bbX \bbq(t) - \tdbbV) +  \bbC \bbq(t)\!=\! \bbV(t)-\bbmu +\bbC\bbq(t). \label{grafn}
\end{align}
Hence, the $j$-th entry of the gradient $\nabla f(\bbq(t))$ does not depend on the full vector $\bbq(t)$, but only local information on its own bus voltage $V_j(t)$ and VAR input $q_j(t)$.
Hence, the GP iteration \eqref{graprojd} can be completely decoupled into local updates, as given by 
\begin{align}
q_j(t+1) =& [1-\alpha(t)] q_j(t) + \alpha(t) \times \nonumber\\
 &\mathbb P_j\! \left[(1-d_jc_j) q_j(t) - d_j (V_j(t) - \mu_j) \right] ~\forall j  \label{locQ}
\end{align}
where $\mathbb P_j$ denotes the projection at bus $j$ to the interval $[\underline q_j, \bbarq_j]$, which is again a local computation. The proposed local VAR control design relying on \eqref{locQ} is essentially equivalent to the centralized GP solver for \eqref{cenqn}. }

\begin{proposition}\label{prop:loc}
Under constant $\bar{\bbV}$, the fixed-point of the iterative update \eqref{graprojd}, or equivalently its local counterpart \eqref{locQ}, will achieve the optimum $\bbq^\star$ to the VAR control problem \eqref{cenqn}.
%if all the eigenvalues $|1 -\lambda_k(\bbD(\bbX+\bbC))| <1$, ~$\forall k$.
\end{proposition}

{\hao {\it Proof:} 
First, existence and uniqueness of the optimum $\bbq^\star$ follows from the strong convexity and none-empty constraint of \eqref{cenqn}. Furthermore, the first-order optimality condition \cite[Prop. 2.1.2]{Ber_NPbook} for \eqref{cenqn} implies that $[\nabla f(\bbq^\star)]^T (\bbq-\bbq^\star) \geq 0$ holds  $\forall \bbq \in [\underline{\bbq}, \bar{\bbq}]$. Thus, the necessary condition for $\bbq^\star$ to be optimum boils down to
\begin{align}
\frac{\partial f(\bbq^\star)} {\partial q_j} 
\left\{ \begin{array}{cl} = 0 & \mathrm{if}~ \underline q_j  \leq q_j^\star \leq \bbarq_j\\
\geq  0 & \mathrm{if}~q_j^\star = \underline q_j \\
\leq  0 & \mathrm{if}~q_j^\star = \bbarq_j
\end{array}~~, \right.
\end{align}
or equivalently, %the optimum stationary point to \eqref{cenqn} would satisfy
\begin{align}
\bbq^\star &= \mathbb P [\bbq^\star - \bbD \nabla f(\bbq^\star)] \nonumber\\
&=\mathbb P [(\bbI - \bbD \bbC) \bbq^\star-\bbD(\bbX\bbq^\star -\tdbbV) ]
\label{locQs}
\end{align}
would hold for any PD diagonal matrix $\bbD$. This is exactly the condition that a fixed point to \eqref{graprojd} would satisfy.  \qed}

%Existence and uniqueness of the optimum $\bbq^\star$ follows from the strong convexity of \eqref{cenqn},  while Proposition \ref{prop:loc} can be proved by checking the first-order optimality condition \cite[Prop. 2.1.2]{Ber_NPbook}, that is
%\begin{align}
%\bbq^\star &= \mathbb P [\bbq^\star - \bbD \nabla f(\bbq^\star)] \nonumber\\
%&=\mathbb P [(\bbI - \bbD \bbC) \bbq^\star-\bbD(\bbX\bbq^\star -\tdbbV) ].
%\label{locQs}
%\end{align}
%%
{\hao
\begin{remark} \label{rmk:feature} (Features of local control.) {\rm 
The GP-based VAR control \eqref{locQ} requires each bus to measure its local voltage magnitude, which can be implemented with minimal hardware updates. The computational requirement is very minimal since \eqref{locQ} involves only scalar operations. As detailed soon, the GP stepsize parameters need to be chosen judiciously to ensure the dynamic stability. However, they can be determined off-line using solely system topology and line admittance information. Without the network information, it is also possible to develop some parameter tuning schemes that adaptively react to potential local oscillation. In addition, even though the analytical results rely on the linearized power flow model, the local update \eqref{locQ} itself can be easily implemented in any three-phase feeders with lossy and coupling effects. Last but not least, the local VAR control can also be updated in an asynchronous fashion with various update rates under minimal centralized coordination. This is of particular interest for the ``plug-and-play'' functionality in microgrid design.}
\end{remark} 

Albeit its simple design and robust features, the local VAR update could suffer from performance degradation due to the surrogate error norm objective \eqref{cenqfn}. As mentioned earlier, this objective would yield the same optimum solution under unlimited VAR resources. For the case of limited but abundant VAR resources, it is expected that the GP-based local update could closely approximate the desired optimum $\bbq^\dag$ of the original problem \eqref{cenq}. Otherwise, it is possible to quantify the exact performance degradation from numerical simulations with given feeder information as shown in Sec. \ref{sec:sim}.}

\section{Stable Local VAR Control} \label{sec:local}

%The two first local Var control methods select a constant $\alpha(t)=1$, under which the VAR update \eqref{locQ} boils down to pure projection, as given by 
%%
%\begin{align}
%q^g_j(t+1) &= \mathbb P_j [(1-d_jc_j) q^g_j(t) - d_j (V_j(t) - \mu_j) ] ~~\forall j.  \label{locQp}
%\end{align}
%%
%This is the simplest form of the gradient projection methods, with the stepsize $d_j>0$. The convergence of the projection based VAR control depends on the local stepsize $d_j$'s, or namely the diagonal matrix $\bbD \succ \bb0$. This is related to the stability analysis since the iterative update \eqref{locQp} can be viewed as a discrete-time system. Therefore, the convergence of \eqref{locQp} with any initial conditions ensures that the associate dynamic system under \eqref{ldfVm4} is stable. 

 %This also relates to the stability of the local VAR control given by \eqref{locQ} if the system is no longer static, which will be discussed in the ensuring section.
 {\ml We will analyze the dynamic stability issue of the proposed GP-based VAR control in this section.}
With highly variable renewable generation and elastic load, the stability of the local VAR control will be crucial to cope with the system dynamics at a very fast time-scale. This is closely related to the static-case analysis on the convergence conditions of \eqref{locQ}. This is equivalent to having the iterative error mismatch  $[\bbq(t)-\bbq^\star] \rightarrow \bb0$ when $t\rightarrow 0$. To this end, diagonally scale all iterates with $\bbD^{-1/2}=\diag(1/\sqrt{d_1},\ldots, 1/\sqrt{d_N})$, and define $\tdbbq(t) := \bbD^{-1/2} \bbq(t)$ and $\tdbbq^\star := \bbD^{-1/2} \bbq^\star$. Since  projection is a nonexpansive mapping, the scaled error norm 
\begin{align}
&\|\tdbbq(t+1) - \tdbbq^\star\| = \big\|\tdbbq(t+1) - [1-\alpha(t)+\alpha(t)]\tdbbq^\star\big\| \nonumber\\
\leq &[1-\alpha(t)] \left\|\tdbbq(t)-\tdbbq^\star\right\|  + \alpha(t) \left\|\tdbbq(t) - \bbD^{1/2}\nabla f(\bbq(t)) - \tdbbq^\star\right\| \nonumber\\
= &\! [1-\alpha(t)]\! \left\|\tdbbq(t)\!-\!\tdbbq^\star\right\|\!+\! \alpha(t)\left\|[\bbI - \bbD^{1/2}(\bbX+\bbC)] (\bbq(t) \!-\! \bbq^\star) \right\|\nonumber\\
 \leq & \left\{1-\alpha(t)+ \alpha(t)\left\|\bbI - \bbD^{1/2}(\bbX+\bbC) \bbD^{1/2} \right\|\!\right\} \| \tdbbq(t) - \tdbbq^\star\| \label{conv}
\end{align}
where the second equality holds by substituting \eqref{locQs}. Denoting matrix $\bbH := \bbD^{1/2}(\bbX+\bbC) \bbD^{1/2}$, one can establish that its $k$-th eigenvalue
\begin{align}
\lambda_k^H:=\lambda_k(\bbH)  = \lambda_k\left(\bbD^{1/2} (\bbX+\bbC)\bbD^{1/2}\right) >0 \label{lamdbak}
\end{align}
where the last inequality holds because $\bbD$, $\bbX$, and $\bbC$ are all PD. A sufficient stability condition is to ensure the non-negative error $\| \tdbbq(t) - \tdbbq^\star\|$ is contracting at every iteration, which leads to the following proposition.
\begin{proposition}\label{prop:stab}
The local update \eqref{locQ} is guaranteed to be stable if $\bbH$'s largest eigenvalue  $\lambda_{\max}^H <2$. 
\end{proposition}

{\hao {\it Proof:} 
By definition the matrix Euclidean norm $\|\bbI - \bbH\|$ equals to the largest singular value of $\bbI-\bbH$. 
Because matrix $\bbH$ is positive definite [cf. \eqref{lamdbak}], $|1-\lambda_k^H|$ is a singular value of $\bbI-\bbH$, implying the scaling coefficient of \eqref{conv}
\begin{align}
1-\alpha(t)+\alpha(t)\|\bbI - \bbH\| \! =\!1-\alpha(t)+\alpha(t)  \max_k|1-\lambda_k^H| ,~\forall t. \nonumber
\end{align}
Assuming $\alpha(t)>\epsilon>0$, having  $|1-\lambda_k^H| <1$ for every $k$ ensures that \eqref{conv} is a contraction mapping for every $t$. Accordingly, the scaled error norm $\| \tdbbq(t) - \tdbbq^\star\|$ will go to 0 in the limit; and similarly for $\|\bbq(t)-\bbq^\star\|$ since $\bbD$ is PD.  
%which leads to that  $|1-\lambda_k^H| <1$ for every $k$ assuming $\alpha(t)>\epsilon>0$. 
Hence, to ensure the stability of \eqref{locQ} the largest eigenvalue of $\bbH$ needs to be less than 2. \qed}

Under the linearized model \eqref{ldfVm4}, Proposition \ref{prop:stab} coincides with the classical dynamic system stability conditions on the eigenvalues of the Jacobian matrix used by \cite{Rob_tps13,JahAli_tps13}. However, our proof based on contracting error norm could handle the projection operator under limited VAR resources.

\subsection{Droop VAR Control}\label{sec:droop}
\begin{figure}[tb]
\centering
\vspace{-3pt}
\includegraphics[width=.66\linewidth,height=1.3in,clip = true, trim = 1.8in 1.88in  1.8in  1.46in]{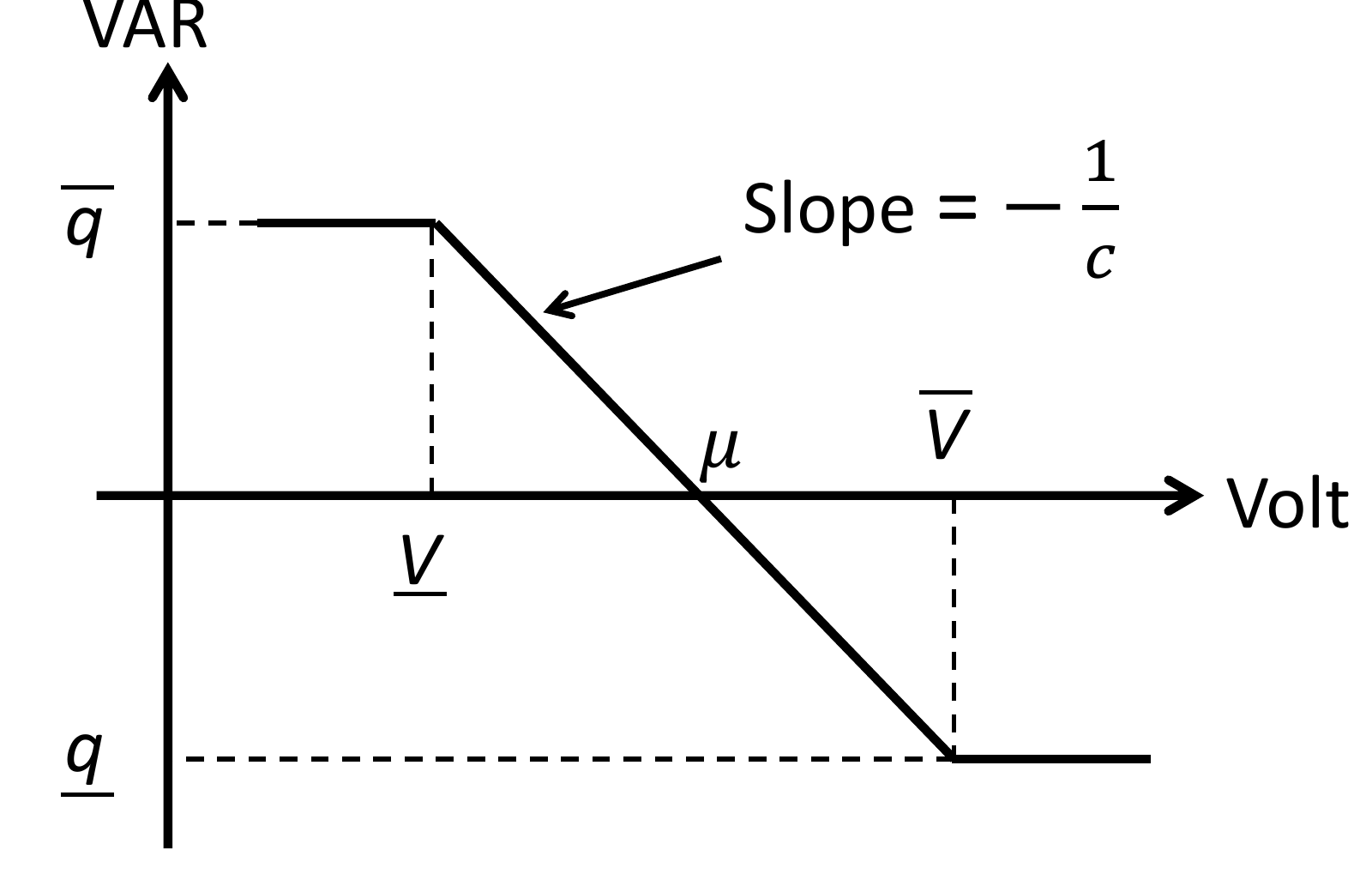}
\caption{Droop VAR control curve.}
\label{fig:droop}
\end{figure}

As advocated by the EPRI smart-inverter initiative \cite{smartinverter} and IEEE 1547.8 standard \cite{ieee1547}, the droop control scheme is to scale the inverter VAR output based on the instantaneous bus voltage mismatch. For given $\bbC$, the GP-based local control \eqref{locQ} actually generalizes this scheme by setting  $d_j = 1/c_j$ and $\alpha(t)=1$, yielding
\begin{align}
q_j(t+1) = \mathbb P_j [- c_j^{-1} (V_j(t) - \mu_j) ], ~~\forall j.  \label{droopQ}
\end{align}
Such a linear droop control design has a negative slope of $-c_j^{-1}$ and no deadband, as illustrated by Fig. \ref{fig:droop}. Based on \eqref{costqg}, the larger $c_j$ is, the more costly it is to provide VAR at bus $j$. Accordingly, the droop slope would be smaller in magnitude to make it less sensitive to the local voltage mismatch. As pointed out by \cite{JahAli_tps13}, this local control scheme is prone to instability, while \cite{FarLow_cdc13} has proven that \eqref{droopQ} is stable {\hao if $(\bbC^{-1} - \bbX)$ is PD}. Proposition \ref{prop:stab} also generalizes this result, since the droop condition $\bbD= \bbC^{-1}$ leads to
\begin{align}
\lambda_{\max}^H %= \lambda_{\max}(\bbC^{-1/2}\bbX\bbC^{-1/2}+\bbI) 
= 1+\lambda_{\max}(\bbC^{-1/2}\bbX\bbC^{-1/2}) < 2, \label{droop_conv}
\end{align}
which is equivalent to {\hao the condition that $(\bbC^{-1}- \bbX)$ is PD because both matrices are PD}. 

If $\bbX$ has very large eigenvalues, the local droop slope needs to be very small and thus it is less sensitive to the voltage mismatch. A smaller droop slope could affect the algorithmic convergence rate. In addition, this would lead to a large $\bbC$ and higher penalty on VAR supply, which is not necessarily desirable. As pointed out by \cite{bz_naps13}, matrix $\bbX$ becomes much worse conditioned if the size of network grows up or the main feeder line gets longer. Hence, it is expected that the droop control strategy would be more likely to be unstable for large-scale distribution systems.

\subsection{Scaled VAR Control}\label{sec:scaled}

To address the stability concerns of droop control, we will consider other designs of $\bbD$ for local VAR control. The Newton's second-order method has been a very popular approach to accelerate the convergence of iterative solvers, by scaling the gradient with the inverse Hessian of the objective; see e.g., \cite[Sec. 1.3]{Ber_NPbook}. However, such a scaling method would be problematic when it comes to constrained problem solvers with a projection operator. It can be easily shown that the fixed-point of the Newton's update generally does not attain the optimum solution of a constrained problem that it intends to solve \cite[Sec. 2.4]{Ber_NPbook}.  Hence, we will set $\alpha(t)=1$ for the iteration \eqref{graprojd}, and scale it using the inverse of the diagonals of Hessian matrix; i.e.,  %\cite[Sec. 2.3]{Ber_NPbook}, as given by
\begin{align}
\bbD = \epsilon~ \bbD_H := \epsilon~[\diag(\bbX + \bbC)]^{-1} \label{scalD}
\end{align}
for some $\epsilon >0$. %Since $\bbD$ is positive and diagonal, we have $\lambda_k(\bbD(\bbX+\bbC)) = \lambda_k(\bbD^{1/2} (\bbX+\bbC)\bbD^{1/2}) >0$ for every $k$. 
The stability condition of Proposition \ref{prop:stab} now becomes
\begin{align}
\epsilon < 2/\lambda_{\max}(\bbD_H^{1/2} (\bbX+\bbC)\bbD_H^{1/2}). \label{scal_conv}
\end{align}
%

%Clearly, design and analysis of a uniform control scheme are also included by setting $\bbD_H = \bbI$. 
{\hao Compared to the droop control, the proposed scaled design can stabilize the system dynamics for any matrix $\bbC$. This offers better flexibility for choosing the VAR supply penalty. In addition, the Hessian based diagonal scaling in \eqref{scalD} helps improve the matrix conditioning and would speed up the convergence rate. More numerical simulations will be given in Sec. \ref{sec:sim} to demonstrate this improvement.}
%As demonstrated by numerical simulations, this scaled VAR control offers better convergence and stability conditions as compared to the droop control. 

\subsection{Delayed VAR Control}\label{sec:delayed}

With a non-unit stepsize $\alpha(t)<1$ in \eqref{locQ}, the most general GP update takes the weighted average between the last iterate $q_j(t)$ and the projection result. This coincides with the delayed droop control method developed in \cite{JahAli_tps13}, in the form of 
\begin{align}
q_j(t+1) \! =\! [1-\alpha(t)]q_j(t)\!+\!\alpha(t)\mathbb P_j [- d_j (V_j(t) - \mu_j) ]  \label{ddQ}
\end{align}
where $d_j$ can be chosen using either the droop or the scaled control design.
The work in \cite{JahAli_tps13} proposes this practical solution to address the instability issues of droop control, along a very general stability condition.  The latter can only be used to check a specific distribution system with all case information given, but does not provide the exact bounds on the stepsize based on the graph based matrices as in Proposition \ref{prop:stab}. 

Convergence of the GP method could also depend on the choice of $\alpha(t)$, which is not reflected by Proposition \ref{prop:stab}. As  in \cite[Sec. 2.3]{Ber_NPbook}, the GP method in the form of \eqref{graprojd} is convergent as long as $\bbD$ is kept constant and $\alpha(t)$ is selected using the limited minimization rule or the Armijo rule. The gist of both rules is to choose $\alpha(t)$ adaptively to ensure sufficient decent in the  objective value at every iteration. Numerical tests performed in \cite{JahAli_tps13} demonstrate that a small and constant choice such as $\alpha(t)= 0.1$, would lead to stable VAR control empirically, even if the corresponding droop control is unstable. This suggests that a small $\alpha(t)$ could contribute to diminish the error norm  in \eqref{conv} as well. 

To better understand the effect of $\alpha(t)$, let us assume there are abundant VAR resources, and thus the projection in \eqref{locQ} is never  active. This is exactly the scenario where instability of the droop control would emerge, as argued at the end of Sec. \ref{sec:droop}. Under this assumption, a closer look at the mismatch error in \eqref{conv} yields
\begin{align}
&\|\tdbbq(t+1) - \tdbbq^\star\| 
 \leq  \big\{1-\alpha(t)\left\|\bbH \right\|\big\} \| \tdbbq(t) - \tdbbq^\star\|. \label{convdelay}
\end{align}
The effective Jacobian now becomes $\alpha(t)\bbH$ for the local control update \eqref{locQ}. Hence, the sufficient stability condition for \eqref{locQ} is updated as $\lambda_{\max}^H <2/\alpha(t)$.  As $\alpha(t)\rightarrow 0$, the delayed VAR control is more likely to be stable.  This analysis corroborates that the stability of droop control can be improved by the delayed design, as pointed out by \cite{JahAli_tps13}.

%%%%%%%%%%%%%%%%%%%%%%%%%%%%%%%%%%%%%%%%%%%%%%%%%%%%%%%%%%%%%%%%%%%%%
 %%
 %%      Section: TC %%
 %%%%%%%%%%%%%%%%%%%%%%%%%%%%%%%%%%%%%%%%%%%%%%%%%%%%%%%%%%%%%%%%%%%%
%%
\section{Numerical Tests} \label{sec:sim}

This section presents numerical test results to demonstrate the effectiveness of the proposed local control methods, for single- and three-phase feeder systems. Both static and dynamic scenarios on the system loading and generation will be considered. 
To better compare various algorithms, the desired voltage  magnitude $\mu_j$ is set to be unit in p.u. at every bus $j$, along with $V_0$ fixed at 1. Each bus is equipped with a certain amount of PV panels, which are able to offer flexible VAR supply to the feeder via effective inverter design. For the dynamic simulation scenarios, the VAR limits $[\underline{\bbq}, \bar{\bbq}]$ are updated at every time slot based on the given inverter ratings and the instantaneous real power generated. 

All numerical tests use the open source simulator OpenDSS \cite{opendss} to solve for the actual power flow, instead of the approximate solution using \eqref{ldfVm3}.  In addition, the actual bus voltage magnitude, instead of the one obtained by the LinDistFlow model, is used for updating the VAR control outputs and numerical performance comparisons.

\subsection{Single-Phase 16-Bus Radial Feeder}

A 12kV radial distribution feeder of 16 buses is first considered; i.e., {\ml the network in Fig. \ref{fig:radial} with $N=15$}. Each line segment has the same impedance of $(0.466+j0.733)\Omega$. For the static case, each bus has a constant load of $(100+j50)$kVA, and abundant VAR resources of $q_j^g \in [-100,~100]$kVA. To include the VAR supply penalty, $c_j$ is chosen to be 0.2 at every bus for the proposed scaled and delayed schemes. The (delayed) droop control will used $c_j=0.5$ based on a linear droop curve of no deadband as in Fig. \ref{fig:droop}, since the voltage limits are set at $[0.95,~1.05]$ and VAR limits at $[-100,~100]$kVA. 

\subsubsection{Static scenario}  
\begin{figure}[tb]
\centering
\includegraphics[width=\linewidth,clip = true, trim = .3in 0.18in  .8in  0.3in]{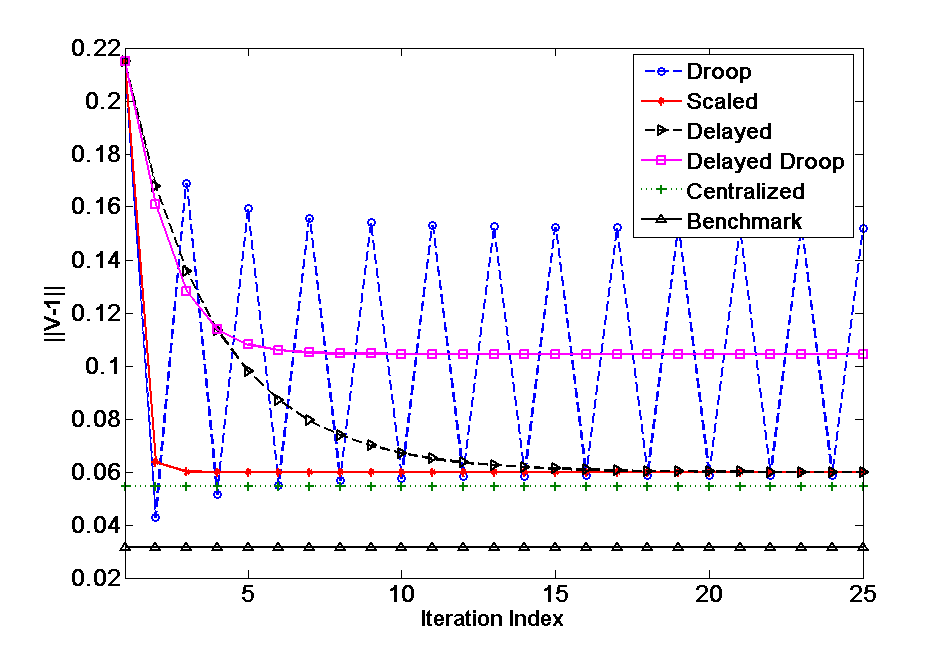}
\caption{Voltage mismatch error versus iteration index for various VAR control methods under the static system setting.}
\label{fig:voltcomp}
%\vspace{-5mm}
\end{figure}

Fig. \ref{fig:voltcomp} plots the iterative voltage mismatch error norm $\|\bbV-\mathbf 1\|$ for the local VAR control methods. The centralized solution corresponds to the optimum solution $\bbq^\star$ to \eqref{cenqn} with $c_j=0.2$. To provide the benchmark performance under the original unweighted objective, the matrix $\bbB$ in \eqref{cenqfn} is substituted by $\bbarlambda \bbI$, where the scalar $\bbarlambda$ is the average of $\bbB$'s eigenvalues. % Both problems are solved with the actual voltage profile at no VAR $\bar{\bbV}$ from AC power flow to correctly capture the system model. 

As for the stepsize choice, $\epsilon$ is chosen to be 0.3 for the scaled and delayed control methods, based on the stability conditions in Proposition \ref{prop:stab}. The stepsize  $\alpha(t)$ is set to be constant 0.3 at every iteration for both delayed schemes.  Choices of stepsize will be discussed soon in more detail. As depicted by Fig. \ref{fig:voltcomp}, the droop control fails to converge to a fixed operating condition, which oscillates between two operating points. This coincides with the earlier discussion that a larger penalty coefficient $\bbC$ due to a steeper linear droop slope would lead to instability. By delaying the droop VAR update with a small $\alpha(t)$, it is possible to stabilize the droop control since the effective Jacobian's eigenvalues will become  smaller. Since the VAR provision penalty is higher for the delayed droop control, its voltage mismatch error will be larger than the scaled and delayed ones. The scaled control method converges very fast to the centralized solution, same for the delayed one.  The best voltage mismatch performance for local control strategies under the surrogate VAR objective is around 0.055, as compared to the benchmark performance around 0.031. This speaks for the competitiveness of the proposed totally local methods, which require minimal coordination in selecting stepsize.

\begin{figure}[tb]
\centering
\includegraphics[width=\linewidth,clip = true, trim = .3in 0.18in .8in  0.3in]{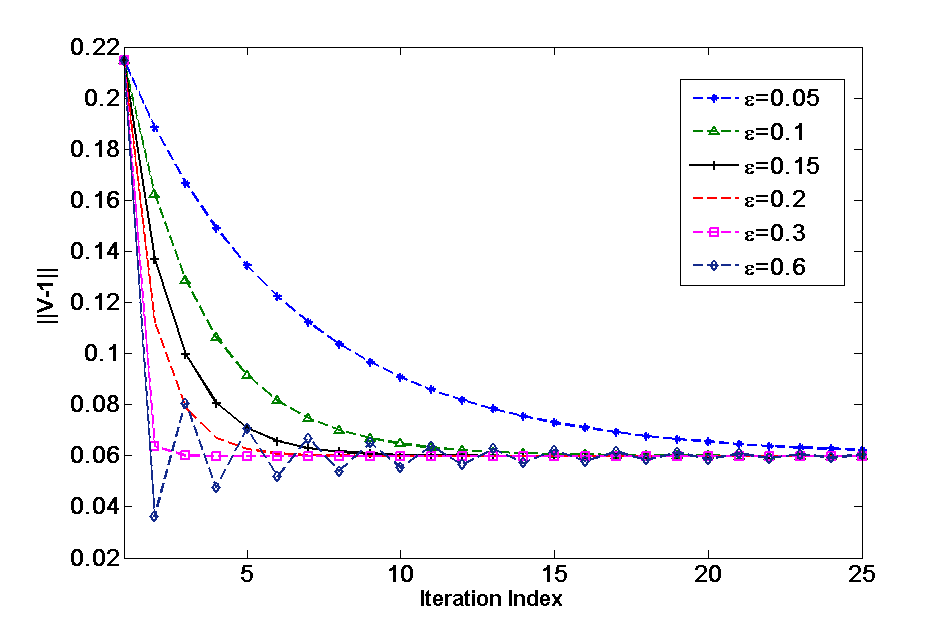}
\caption{Iterative voltage mismatch error performance for the scaled VAR control method with different $\epsilon$ values.}
\label{fig:scaledeps}
\end{figure}

\begin{figure}[tb]
\centering
\includegraphics[width=\linewidth,clip = true, trim = .3in 0.18in  .8in  0.3in]{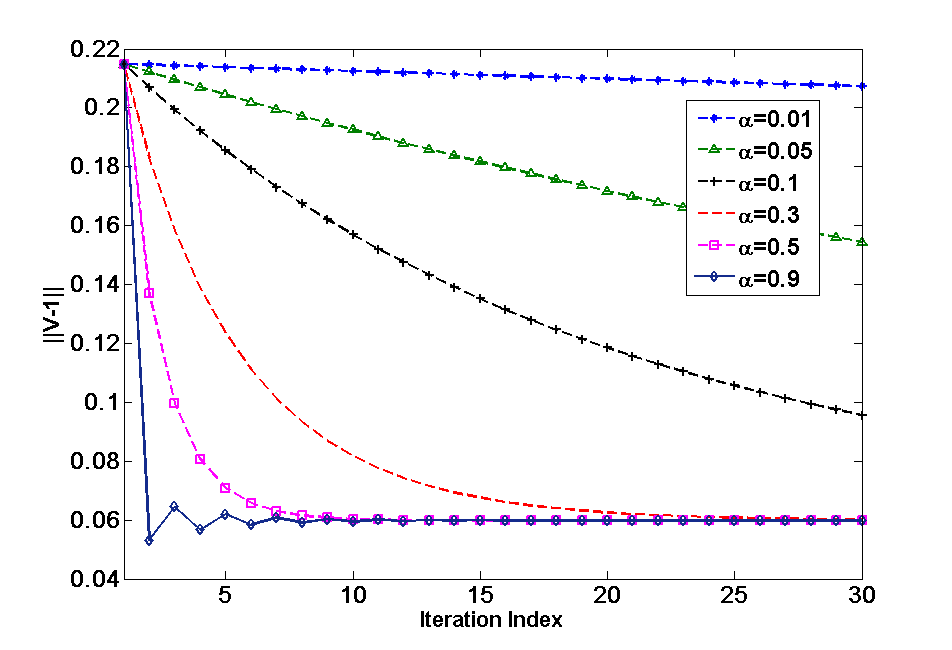}
\caption{Iterative voltage mismatch error performance for the delayed VAR control method with different constant $\alpha$ values.}
\label{fig:delaysalp}
%\vspace{-5mm}
\end{figure}

The effect of stepsize choices on VAR updates is further investigated, which shows a trade-off between the stability and convergence rate. Under the same static setting,  Fig. \ref{fig:scaledeps} plots the voltage mismatch performance for the scaled control method using various $\epsilon$ values. Clearly, the larger $\epsilon$ is, the faster the update converges. However, this could also lead to potential oscillation in the error performance, demonstrating potential instability under fast dynamics. The upper bound of $\epsilon$ is calculated to be around 0.63 according to Proposition \ref{prop:stab}, which coincides with the observation that the error oscillation would happen when $\epsilon = 0.6$. The stepsize of $\epsilon=0.3$ seems to be a very good choice for the 16-bus feeder, while the analysis on the conditioning number of the Jacobian matrix also suggests that this would lead to a very good convergence rate. Hence, the stability results of Sec. \ref{sec:local} are very helpful to select a good $\epsilon$ value if the full feeder information is available. Otherwise, it is also possible to select the stepsize on-the-fly, by adaptively reducing the $\epsilon$ value for each inverter based on its local voltage oscillation intensity. %This further supports the minimal coordination overhead required by the proposed local control schemes. 

Similar analysis has been performed for the delayed control scheme with different constant $\alpha(t)$ values, as plotted in Fig. \ref{fig:delaysalp}. All the error curves are based on fixing $\epsilon=0.3$ as in Fig. \ref{fig:scaledeps}. The value of $\alpha(t)$ is set to be constant for every iteration, varying from 0.01 to 0.9. Under this setting, Fig. \ref{fig:delaysalp} shows that a smaller $\alpha(t)$ would lead to slower convergence rate. It is true that under large $\epsilon$ values, a small $\alpha(t)$ would be very helpful to stabilize the system. Nonetheless, it does not seem that the delayed scheme with the best $\alpha(t)$ would converge faster than the scaled one. Therefore, the numerical tests suggest that the proposed scaled VAR control offers very competitive convergence rate by tuning up the stepsize, which should be advocated for practical implementations. 

{\hao The average algorithmic run time  has been collected using MATLAB$^\circledR$ R2014a software on a typical Windows 8 computer with a 2.6GHz CPU. The computational time for all local (droop, scaled, or delayed) control schemes is around 1-2 microseconds per node per iteration. This numerically corroborate the fast computation feature in Remark \ref{rmk:feature}.}

\begin{figure}[tb]
\centering
\includegraphics[width=\linewidth,clip = true, trim = .38in 0.05in  .6in  0.3in]{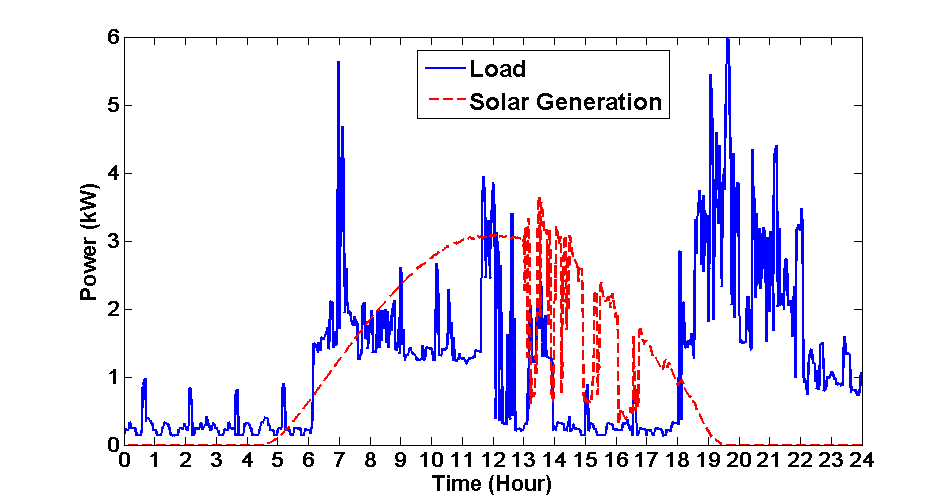}
\caption{Sample daily residential load and solar PV generation profiles.}
\label{fig:load}
%\vspace{-2mm}
\end{figure}

\bigskip
\subsubsection{Dynamic tests} 
Numerical tests under dynamic system operating conditions have also been performed. Dynamic scenarios include PV generation drops due to cloud coverage, and suddenly increasing load during large appliance start-up. To model these, a daily residential load profile and a solar PV generation profile are used, both at minute resolution as shown by Fig. \ref{fig:load}. They are taken from a real dataset in \cite{M_Lich}, which was collected at a residential home in the central Illinois region during a particular day in Summer 2010. {\hao Note that the profile data was collected on a Friday, when house residents are more likely to stay late in the night. This explains why the midnight load is observed to be slightly higher than the house base load in the early morning or afternoon.} The installed solar panels have 3kW peak capacity. The dynamic tests construct the load at every bus to consist of 18 residential homes. Each home has the same PV generation profile as in Fig. \ref{fig:load}, where the inverter apparent power limit is $5\%$ higher than the peak capacity of 3kW. 

\begin{figure}[tb]
\centering
\includegraphics[width=\linewidth,clip = true, trim = .3in 0.0in  .6in  0.25in]{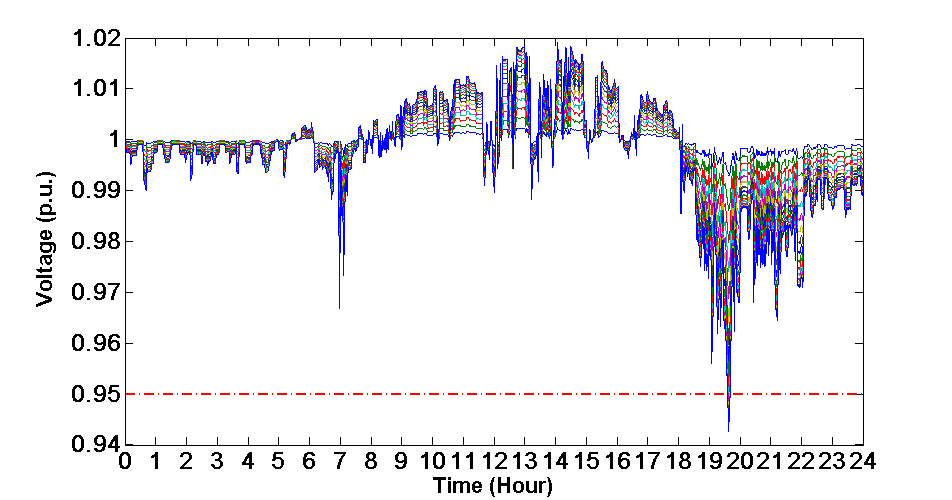}
\caption{Daily bus voltage profile under no VAR support.}
\label{fig:novar}
%\vspace{-5mm}
\end{figure}

\begin{figure}[tb]
{\centering
\includegraphics[width=\linewidth,clip = true, trim = .3in 0.1in  .6in  0.3in]{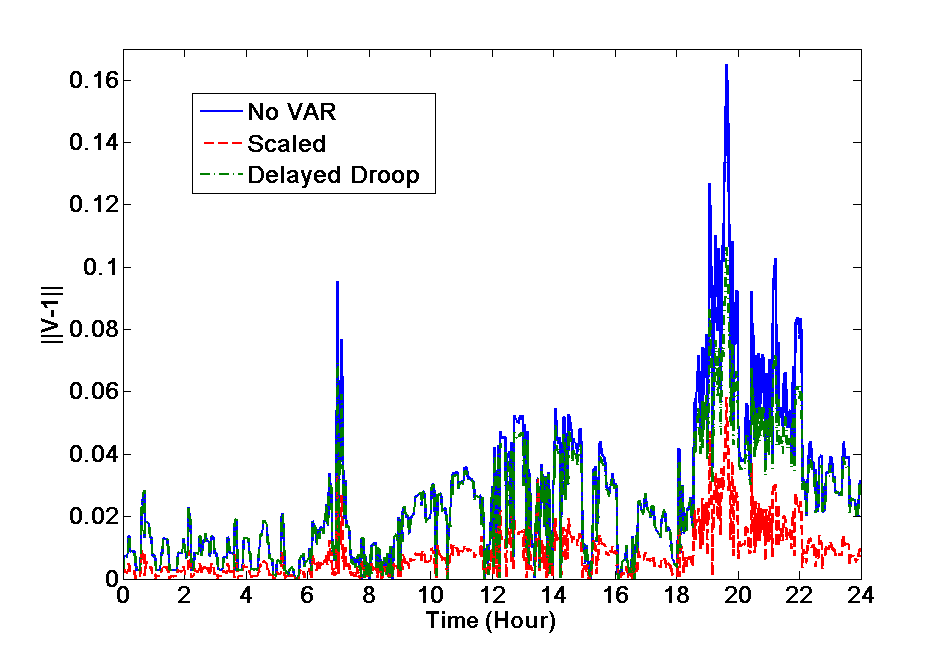}
\caption{Daily voltage mismatch error at every minute for the 16-bus case.}
\label{fig:normvoltage}
}
\end{figure}

Fig. \ref{fig:novar} illustrates the daily feeder voltage profile without any VAR support. {\ml Peak voltage} is observed to happen during noon time when the load is minimal and the solar generation is at its peak. Meanwhile, under-voltage violation (below 0.95 p.u.) has been experienced at the end of feeder, during the evening time with increasing power demand and decreasing PV generation. {\hao Because of its stability issue, the droop method is not suitable for dynamic implementation. In addition, the delayed method has been omitted since it achieves the same steady-state performance as the scaled method but at a slower convergence rate.} Hence, the proposed scaled VAR control scheme is implemented here for tackling the dynamic under- and over-voltage issues, along with the delayed droop method. The parameter and stepsize settings follow from the static tests, while the droop slope is time-varying based on the instantaneous VAR limits computed from the PV generation at every minute. Both local control schemes update every 5 seconds while the load and PV generation stay constant within a minute. The voltage mismatch comparison is illustrated in Fig. \ref{fig:normvoltage}, which shows that additional VAR control outperforms the case with no VAR support. Moreover, because of the constant $c_j=0.2$ setting, the scaled control method is more effective in maintaining a flat voltage profile than the delayed droop design, especially for the evening hours from 18:00 to 22:00. %The delayed control scheme has similar performance as the scaled one, and hence its mismatch error curve is not plotted here. %Hence, its error curve is not plotted here.  

\begin{figure}[tb]
{\hao
\centering
\includegraphics[width=\linewidth,clip = true, trim = .3in 0.1in  .6in  0.3in]{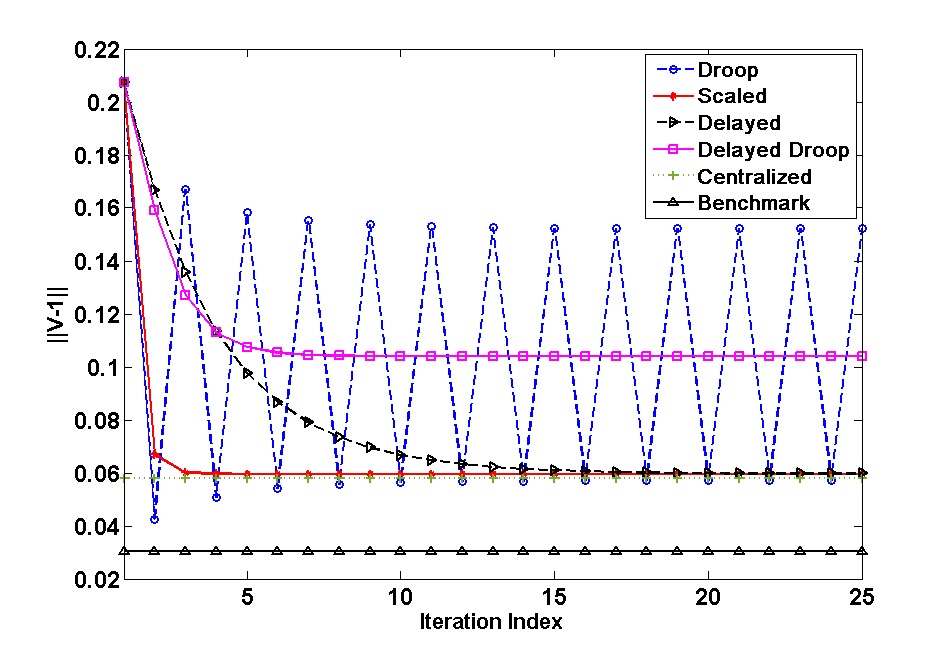}
\caption{Voltage mismatch error versus iteration index for various VAR control methods for the modified 16-bus meshed network.}
%\vspace{-5mm}
\label{fig:compmeshed}
}
\end{figure}

{\hao
\subsection{Modified 16-Bus Meshed Network}
To corroborate Remark \ref{rmk:meshed} on the applicability to meshed networks, we have modified the 16-bus radial feeder by adding two additional lines, one connecting nodes 12-14 and one connecting 13-15. All other system settings follow from the static scenario in Sec. \ref{sec:sim}-A. Fig. \ref{fig:compmeshed} plots the voltage mismatch error per iteration for all control schemes, which shows similar comparison results as in Fig. \ref{fig:voltcomp}. Because the penalty coefficient $c_j = 0.5$ is not sufficiently large, the droop control has shown to be oscillating as well, which has been stabilized using the delayed strategy. The scaled control scheme converges very fast to the centralized optimum solution to the weighted problem \eqref{cenqn}, whose performance slightly degrades from the benchmark solution to the original unweighted problem. This test verifies the applicability of the proposed framework to networks of general topology. 
}

\subsection{IEEE 123-Bus Test Feeder}
\begin{figure}[tb]
\centering
\includegraphics[width=.75\linewidth,height=2in]{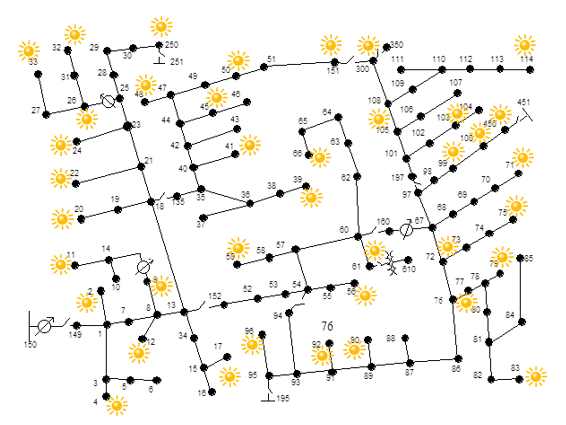}
\caption{IEEE 123-bus test feeder with solar PV panel locations.}
\label{fig:ieee123}
%\vspace{-5mm}
\end{figure}
\begin{figure}[tb]
\centering
\includegraphics[width=.86\linewidth,height=1.68in]{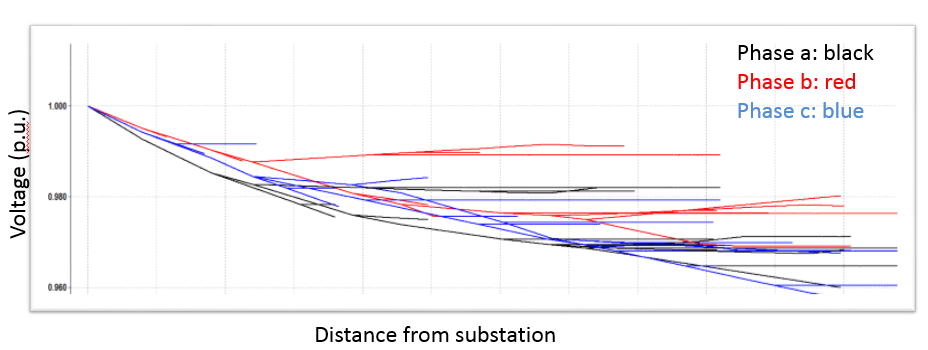}
\caption{Three-phase 123-bus feeder voltage profile.}
\label{unvoltage}
\end{figure}

To demonstrate the effectiveness of the proposed schemes for practical three-phase feeders, we have also implemented them on the IEEE 123-bus test case \cite{ieee123}. Fig. \ref{fig:ieee123} shows the one-line diagram for this distribution feeder case. In order to show the effects of inverter VAR control, the four three-phase voltage regulators are taken out of the 123-bus feeder system. In addition, the case load information is used to determine the number of residential homes connected to each load bus, while each home's load demand and PV generation profile are same as Fig. \ref{fig:load}. %The residential load and PV generation profiles follow from Fig. \ref{fig:load}. 
Locations of the load buses that are equipped with solar panels are shown in Fig. \ref{fig:ieee123}. The feeder voltage profile for this three-phase system under no VAR support is plotted in Fig. \ref{unvoltage} for one time instance, demonstrating that the system is unbalanced. 

\begin{figure}[tb]
{\hao
\centering
\includegraphics[width=\linewidth,height= 2.5in,clip = true, trim = .3in 0.1in  .6in  0.3in]{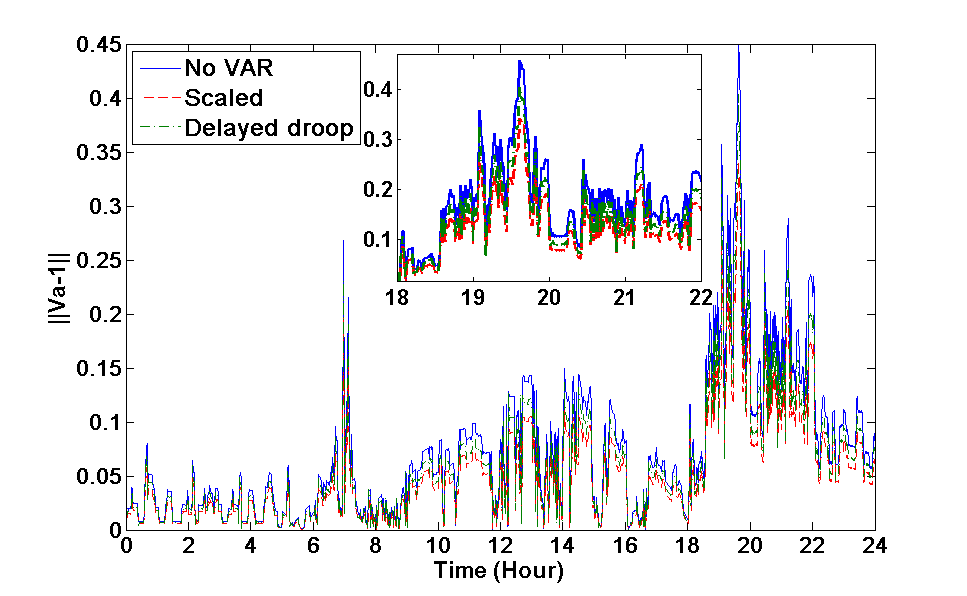}
\centerline{(a)}
\includegraphics[width=\linewidth,clip = true, trim = .3in 0.1in  .6in  0.3in]{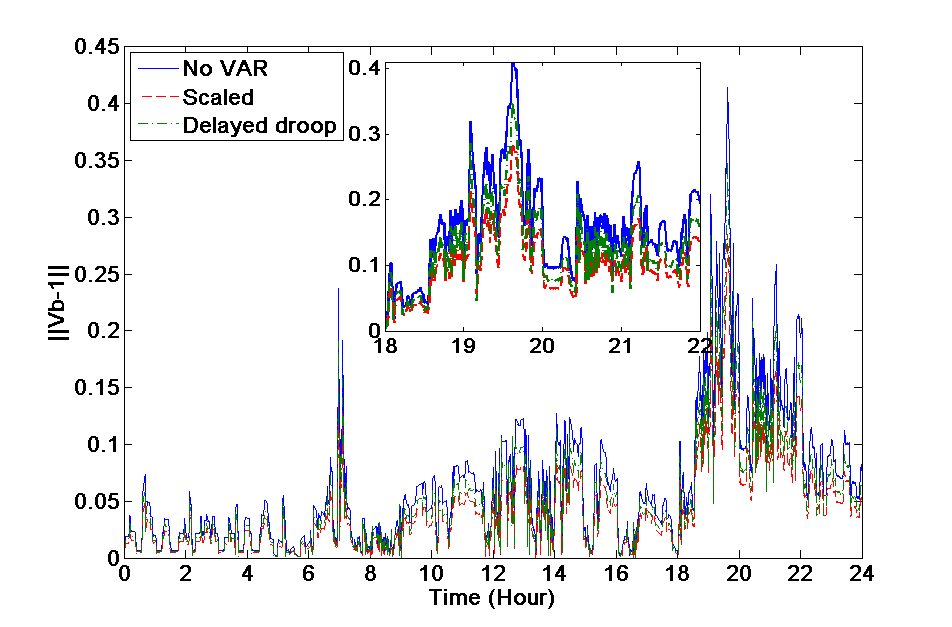}
\centerline{(b)}
\includegraphics[width=\linewidth,clip = true, trim = .3in 0.1in  .6in  0.3in]{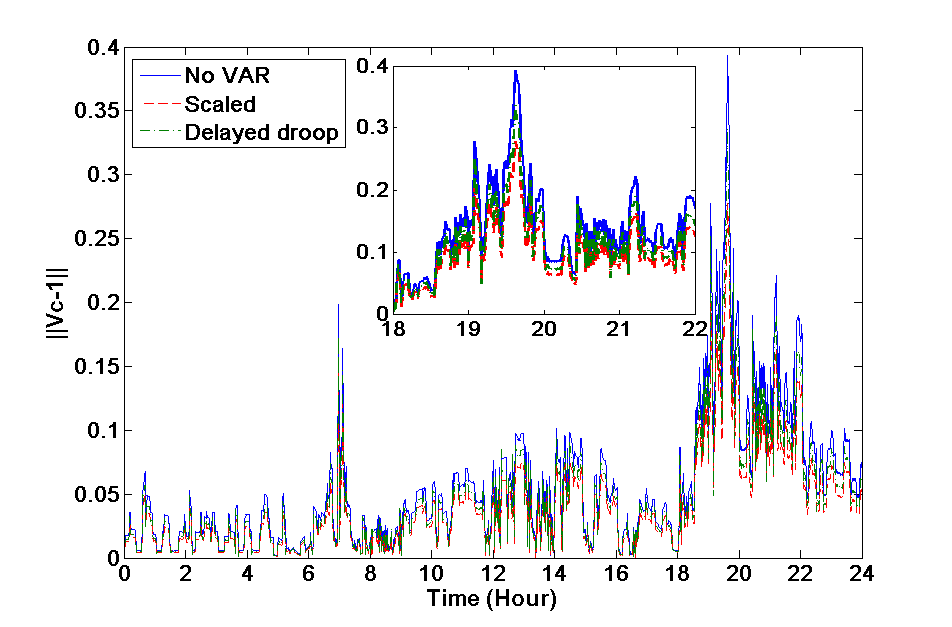}
\centerline{(c)}
\caption{Daily voltage mismatch error at (a) phase a; (b) phase b; and (c) phase c for the 123-bus case.}
\label{fig:misvoltage300}}
\end{figure}

{\hao Figs. \ref{fig:misvoltage300}(a)-(c) plot the daily voltage mismatch error at all three phases for the 123-bus case. Similar to Fig. \ref{fig:normvoltage} for the 16-bus case, only the scaled and delayed droop control methods are compared to the case of no VAR support. For the scaled control, $c_j$ is set to be constant at 0.01 while the delayed droop slope is again time-varying based on the instantaneous inverter VAR limits. Both local VAR control methods improve the voltage support performance over the no VAR support scenario. Because of its constant minimal VAR penalty, the scaled VAR control scheme slightly outperforms the delayed droop method. However, the difference between voltage mismatch error among all three scenarios is less significant compared to the 16-bus system in Fig. \ref{fig:normvoltage}. This less significant voltage regulation performance is because the 123-bus case has a lower level of PV penetration (30\%) compared to the 16-bus case (100\%). Accordingly, the inverter VAR output would be less effective in attaining the ideal constant voltage profile. Nonetheless, the scaled VAR control becomes more effective during the evening hours at higher voltage violation. As shown by the zoom-in view in Figs. \ref{fig:misvoltage300}(a)-(c), the scaled scheme almost reduces the voltage mismatch error by half from the baseline case of no VAR support. }

%%%%%%%%%%%%%%%%%%%%%%%%%%%%%%%%%%%%%%%%%%%%%%%%%%%%%%%%%%%%%%%%%%%%%
 %%
 %%      Section: conclusions %%
 %%%%%%%%%%%%%%%%%%%%%%%%%%%%%%%%%%%%%%%%%%%%%%%%%%%%%%%%%%%%%%%%%%%%
%%

\section{Conclusions and Future Work} \label{sec:con}

This paper presents a general framework for developing local VAR control methods with high penetration of distributed VAR resources. By linearizing the distribution  network power flow model, the VAR control problem is cast as one to minimize the voltage mismatch error. Using the graph matrix representation, we formulate a weighted error minimization problem under box constraints that represent VAR limits at every bus. The gradient-projection (GP) method is evoked for solving this constrained problem, which naturally decouples into local VAR updates requiring only the instantaneous bus voltage magnitude information. The GP-based VAR control framework generalizes existing droop and delayed droop control methods, while allows for stability analysis to tune up the parameters based on the network Bbus matrix. Numerical tests have been performed on single- and three-phase systems using the exact ac power flow model, which corroborate the analytical results on the performance guarantees of proposed methods for realistic system implementations.  

The future research plan for this work includes to investigate
the impact of potential asynchronous control updates among different buses, due to the potential lack of coordination \cite{hznl_icassp16}. We are also actively {\ml investigating the interactions between inverter VAR control and other voltage regulation devices}, as well as pursuing a distributed VAR control framework which has the potential to achieve the original unweighted error objective.

%%%%%%%%%%%%%%%%%%%%%%%%%%%%%%%%%%%%%%%%%%%%%%%%%%%%%%%%%%%%%%%%%%%%%%
%                                                                    %
%       Bibliography %
%
%%%%%%%%%%%%%%%%%%%%%%%%%%%%%%%%%%%%%%%%%%%%%%%%%%%%%%%%%%%%%%%%%%%%%%
\begin{comment}

\end{comment}

\bibliographystyle{IEEEtran}

\itemsep2pt
\bibliography{localvvc,hzpub}

%%\newpage
%\section*{Biographies}
%
%\begin{IEEEbiography}[{\includegraphics[width=1in,
%height=1.25in,keepaspectratio]{haozhu.eps}}]{Hao Zhu}
%(M'12) is currently an Assistant Professor of ECE at the University of Illinois, Urbana-Champaign (UIUC). She received her B.S. from Tsinghua University, Beijing, China in 2006, and the M.Sc. and Ph.D. degrees from the University of Minnesota, Minneapolis in 2009 and 2012, respectively. She has been working as a postdoctoral research associate on power system validation with the Information Trust Institute at UIUC from 2012 to 2013. Her current research interests include power system monitoring and operations, dynamics and stability, distribution systems, and energy data analytics.
%\end{IEEEbiography}

\end{document}